\documentclass[11pt]{article}
\usepackage[T1]{fontenc}
\usepackage{graphicx}
\usepackage{grffile}
\usepackage{longtable}
\usepackage{wrapfig}
\usepackage{rotating}
\usepackage[normalem]{ulem}
\usepackage{amsmath}
\usepackage{caption}
\usepackage{subcaption}
\usepackage{textcomp}
\usepackage{amssymb}
\usepackage{capt-of}
\usepackage{hyperref}
\usepackage{stmaryrd}
\usepackage{algorithm}
\usepackage{xcolor}
\usepackage{algorithmicx}
\usepackage{algpseudocode}
\author{Tom Gustafsson\footnote{Corresponding author. \texttt{firstname.lastname@aalto.fi}}}
\date{\today}
\hypersetup{
  colorlinks=true,
  linkcolor=blue,
  filecolor=blue,
  citecolor=blue,
  urlcolor=blue,
}
\title{A simple technique for unstructured mesh generation via adaptive finite elements}

\begin{document}

\newcommand*\DNA{\textsc{dna}}

\newcommand*\Let[2]{\State #1 $\gets$ #2}
\algrenewcommand\algorithmicrequire{\textbf{Precondition:}}
\algrenewcommand\algorithmicensure{\textbf{Postcondition:}}

\maketitle

\begin{abstract}
  This work describes a concise algorithm for the generation of triangular
  meshes with the help of standard adaptive finite element methods.  We
  demonstrate that a generic adaptive finite element solver can be repurposed
  into a triangular mesh generator if a robust mesh smoothing algorithm is
  applied between the mesh refinement steps.  We present an implementation of
  the mesh generator and demonstrate the resulting meshes via examples.
\end{abstract}

\section{Introduction}
\label{sec:orge4667b0}

Many numerical methods for partial differential equations (PDE's), such as the
finite element method (FEM) and the finite volume method (FVM), are based on
splitting the domain of the solution into primitive shapes such as triangles or
tetrahedra.  The collection of the primitive shapes, i.e.~the computational
mesh, is used to define the discretisation, e.g., in the FEM, the shape
functions are polynomial in each mesh element, and in the FVM, the discrete
fluxes are defined over the cell edges or faces.

This article describes a simple approach for the triangulation of
two-dimensional polygonal domains.  The process can be summarised as follows:
\begin{enumerate}
\item Find a constrained Delaunay triangulation (CDT) of the polygonal domain
      using the corner points as input vertices and the edges
      as constraints.
\item Solve the Poisson equation with the given triangulation
      and the FEM.
\item Split
      the triangles with the largest error indicator
      using adaptive mesh refinement techniques.
\item Apply centroidal patch tesselation (CPT) smoothing to the resulting
  triangulation.
\item Go to step 2.
\end{enumerate}
It is noteworthy that the steps 2, 3 and 5 correspond exactly to what is done in
any implementation of the standard adaptive FEM;
cf.~Verf\"{u}rth~\cite{Verf_rth_2013} who calls it the \emph{adaptive process}.

The goal of this work is to demonstrate that if the mesh smoothing algorithm of
step 4 is chosen properly, the adaptive process tends to produce reasonable
meshes even if the initial mesh is of low quality.  Thus, we demonstrate that
the adaptive process---together with an implementation of the CDT and a robust
mesh smoothing algorithm---can act as a simple triangular mesh generator.

\section{Prior work}
\label{sec:org7798f6d}

Two popular techniques for generating unstructured meshes are those based on the
\emph{advancing front technique}~\cite{L_hner_1988} or the \emph{Delaunay mesh
refinement}~\cite{Chew_1989, Ruppert_1995, Shewchuk_2002}.  In addition, there
exist several less known techniques such as \emph{quadtree
meshing}~\cite{Yerry_1983}, \emph{bubble packing}~\cite{Shimada_1995}, and
hybrid techniques combining some of the above~\cite{mavriplis1995advancing}.

Some existing techniques bear similarity to ours.  For example,
Bossen--Heckbert~\cite{bossen1996pliant} start with a CDT and improve it by
relocating the nodes.  Instead of randomly picking nodes for relocation, we use
a finite element error indicator that guides the refinement.  Instead of doing
local modifications, we split simultaneously all triangles that have their error
indicators above a predefined threshold.

Persson--Strang~\cite{persson2004simple} describe another technique based on
iterative relocation of the nodes.  An initial mesh is given by a structured
background mesh which is then relaxed by interpreting the edges as a
precompressed truss structure.  The structure is forced inside a given domain by
expressing the boundary using signed distance functions and interpreting the
signed distance as an external load acting on the truss.  In contrast to the
present approach, the geometry description is implicit, i.e.~the boundary is
defined as the zero set of a user-given distance function.

We do not expect our technique to surpass the existing techniques in the quality
of the resulting meshes or in the computational efficiency.  However, the
algorithm can be easier to understand for those with a background in the finite
element method and, hence, it may be a viable candidate for supplementing
adaptive finite element solvers with basic mesh generation capabilities.

\section{Components of the mesh generator}
\label{sec:components}

The input to
our mesh generator is a sequence of \(N\) corner points
$$\mathcal{C} = (\mathcal{C}_1, \mathcal{C}_2, \dots, \mathcal{C}_N), \quad \mathcal{C}_j \in \mathbb{R}^2, \quad j = 1,\dots,N,$$
that form a polygon when connected by the edges
$$(\mathcal{C}_1, \mathcal{C}_2), (\mathcal{C}_2,\mathcal{C}_3),
\dots, (\mathcal{C}_{N-1}, \mathcal{C}_N), (\mathcal{C}_N,\mathcal{C}_1).$$
We do not allow self-intersecting polygons although
the algorithm generalises to polygons
with polygonal holes.
The corresponding domain is denoted
by \(\Omega_{\mathcal{C}} \subset \mathbb{R}^2\).

\subsection{Constrained Delaunay triangulation}
\label{sec:cdt}

A \emph{triangulation} of \(\Omega_{\mathcal{C}}\) is a
collection of nonoverlapping nondegenerate triangles whose union is exactly
\(\Omega_{\mathcal{C}}\).  Our initial triangulation \(\mathcal{T}_0\) is a
\emph{constrained Delaunay triangulation} (CDT) of the input vertices
\(\mathcal{C}\) with the edges \((\mathcal{C}_1, \mathcal{C}_2)\),
\((\mathcal{C}_2,\mathcal{C}_3)\), \(\dots\), \((\mathcal{C}_{N-1},
\mathcal{C}_N)\), \((\mathcal{C}_N, \mathcal{C}_1)\) constrained to be a part of
the resulting triangulation and the triangles outside the polygon removed;
cf.~Chew~\cite{Chew_1987} for the exact definition of a CDT and an algorithm for
its construction.

\begin{figure}[htbp]
  \centering
  \includegraphics[width=0.24\textwidth]{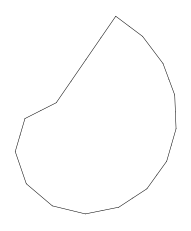}
  \includegraphics[width=0.24\textwidth]{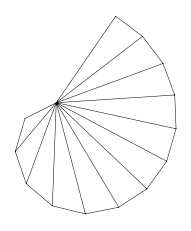}
\caption{A spiral-shaped boundary approximated by linear segments and the
  corresponding CDT with the triangles outside of the polygon removed.  An
  example from the documentation of the Triangle mesh
  generator~\cite{shewchuk1996triangle}.}
\label{fig:cdt}
\end{figure}

An example initial triangulation of a polygon with a spiral-shaped boundary is
given in Figure~\ref{fig:cdt}. It is obvious that the CDT is not always a high
quality computational mesh due to the presence of arbitrarily small angles.
Thus, we seek to improve the initial triangulation by iteratively adding new
triangles, and smoothing the mesh.  Note that the remaining steps do not assume
the use of CDT as an initial triangulation---any triangulation with the
prescribed edges will suffice.

\subsection{Solving the Poisson equation}
\label{sec:poisson}

In order to decide on the placement of the new vertices and triangles, we solve
the Poisson equation\footnote{The choice of the Poisson equation is motivated by
the following heuristic observation: a \emph{quality mesh} is often synonymous
with a \emph{good mesh for the finite element solution of the Poisson
equation}.} using the FEM and evaluate the corresponding a~posteriori error
estimator.  The triangles that have the highest values of the error estimator
are refined, i.e.~split into smaller triangles.

The problem reads: find $u : \Omega_{\mathcal{C}} \rightarrow \mathbb{R}$ satisfying
\begin{alignat}{2}
-\Delta u &= 1 \quad && \text{in $\Omega_{\mathcal{C}}$,} \\
u &= 0 \quad && \text{on $\partial \Omega_{\mathcal{C}}$.}
\end{alignat}
The finite element method is used to numerically solve the
weak formulation: find \(u \in V\) such that
\begin{equation}
   \label{eq:weakform}
   \int_{\Omega_{\mathcal{C}}} \nabla u \cdot \nabla v \,\mathrm{d}x = \int_{\Omega_{\mathcal{C}}} v\,\mathrm{d}x \quad \forall v \in V,
\end{equation}
where
\(w \in V\) if \(w |_{\partial \Omega_{\mathcal{C}}} = 0\) and
$
   \int_{\Omega_{\mathcal{C}}} (\nabla w)^2 \,\mathrm{d}x < \infty.
$

We denote the \(k\)th triangulation of the
domain \(\Omega_{\mathcal{C}}\) by \(\mathcal{T}_k\), \(k=0,1,\dots\), and
use the piecewise linear polynomial space
$$V_h^k = \{ v \in V : v|_T \in P_1(T)~\forall T \in \mathcal{T}_k \},$$
where $P_1(T)$ denotes the set of linear polynomials over $T$.
The finite element method corresponding to the \(k\)th iteration reads:
find \(u_h^k \in V_h^k\) such that
\begin{equation}
   \label{eq:discweakform}
   \int_{\Omega_{\mathcal{C}}} \nabla u_h^k \cdot \nabla v_h \,\mathrm{d}x = \int_{\Omega_{\mathcal{C}}} v_h\,\mathrm{d}x \quad \forall v_h \in V_h^k.
\end{equation}
The local a posteriori error estimator
reads
\begin{equation}
        \eta_T(u_h^k) = \sqrt{h_T^2 A_T^2 + \frac12 h_T \int_{\partial T \setminus \partial \Omega_{\mathcal{C}}} (\llbracket \nabla u_h^k \cdot \boldsymbol{n} \rrbracket)^2 \,\mathrm{d}s}, \quad T \in \mathcal{T}_k,
\end{equation}
where $A_T$ is the area of the triangle $T$ and $h_T$ is the length of its longest edge, $\llbracket w \rrbracket |_{\partial T \setminus \partial \Omega_{\mathcal{C}}}$ denotes the jump in the values of
$w$ over $\partial T \setminus \partial \Omega_{\mathcal{C}}$, and $\boldsymbol{n}$ is a unit normal vector to
$\partial T$.
The error estimator $\eta_T$ is evaluated for each triangle
after solving \eqref{eq:discweakform}.
Finally, a triangle $T \in \mathcal{T}_k$ is marked for refinement if
\begin{equation}
  \label{eq:adaptivetheta}
   \eta_T > \theta \max_{T^\prime \in \mathcal{T}_k} \eta_{T^\prime},
\end{equation}
where $0 < \theta < 1$ is a parameter controlling the amount
of elements to split during each iteration. \cite{Verf_rth_2013}

\subsection{Red-green-blue refinement}
\label{sec:rgb}

\begin{figure}[htbp]
\centering
\includegraphics[width=0.24\textwidth]{./image_5.png}
\includegraphics[width=0.24\textwidth]{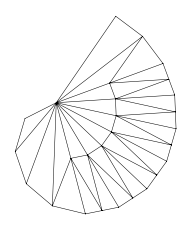}
\caption{(Left.) The initial triangulation. (Right.) The resulting triangulation
  after a solve of \eqref{eq:discweakform} and an adaptive RGB refinement.}
\label{fig:firstrgb}
\end{figure}

The triangles marked for refinement by the rule \eqref{eq:adaptivetheta} are
split into four.  In order to keep the rest of the triangulation conformal,
i.e.~to not have nodes in the middle of an edge, the neighboring triangles
are split into two or three by the so-called red-green-blue (RGB) refinement;
cf.~Carstensen~\cite{carstensen2004adaptive}.  Using RGB refinement to the
example of Figure~\ref{fig:cdt} is depicted in Figure~\ref{fig:firstrgb}.

\subsection{Centroidal patch triangulation smoothing}
\label{sec:cpt}

\begin{figure}[htbp]
\centering
\includegraphics[width=0.24\textwidth]{./image_6.png}
\includegraphics[width=0.24\textwidth]{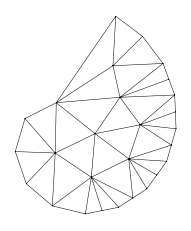}
\caption{(Left.) The adaptively refined triangulation. (Right.) The resulting
  triangulation after smoothing and edge flipping.}
\label{fig:firstsmooth}
\end{figure}

We use a mesh smoothing approach introduced by Chen--Holst~\cite{Chen_2011} who
refer to the algorithm as \emph{centroidal patch triangulation} (CPT) smoothing.
The idea is to repeatedly move the interior vertices to the area-weighted
averages of the barycentres of the surrounding triangles.  The CPT smoothing is
combined with an edge flipping algorithm, also described in
Chen--Holst~\cite{Chen_2011}, to improve the quality of the resulting
triangulation.  The mesh smoother is applied to the spiral-shaped domain example 
in Figure~\ref{fig:firstsmooth}.

\section{The mesh generation algorithm}
\label{sec:orgff9b6c1}

\begin{figure}[htbp]
  \centering
  \includegraphics[width=0.24\textwidth]{./image_19.png}
  \includegraphics[width=0.24\textwidth]{./image_5.png}
  \includegraphics[width=0.24\textwidth]{./image_6.png}
  \includegraphics[width=0.24\textwidth]{./image_7.png}\\
  \includegraphics[width=0.24\textwidth]{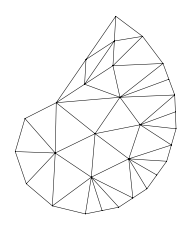}
  \includegraphics[width=0.24\textwidth]{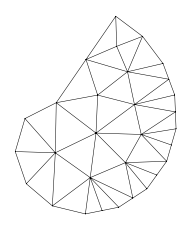}
  \includegraphics[width=0.24\textwidth]{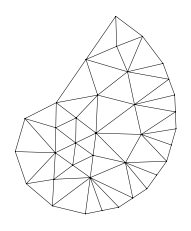}
  \includegraphics[width=0.24\textwidth]{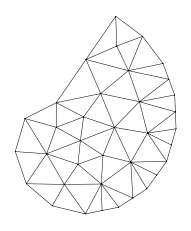}\\
  \includegraphics[width=0.24\textwidth]{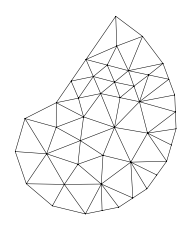}
  \includegraphics[width=0.24\textwidth]{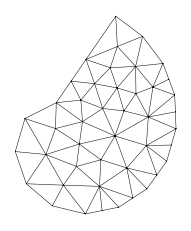}
  \includegraphics[width=0.24\textwidth]{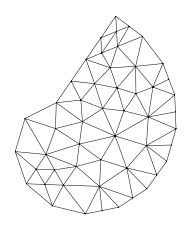}
  \includegraphics[width=0.24\textwidth]{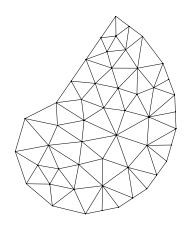}\\
  \includegraphics[width=0.24\textwidth]{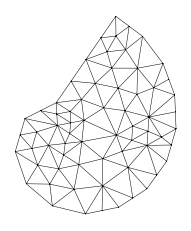}
  \includegraphics[width=0.24\textwidth]{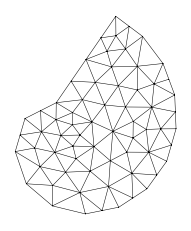}
  \hspace{0.48\textwidth}
  \caption{The entire mesh generation process for the spiral-shaped domain
    example from left-to-right, top-to-bottom.}
\label{fig:spiralexample}
\end{figure}

In previous sections we presented an overview of all the components of the mesh
generation algorithm.  The resulting mesh generator is now summarised in
Algorithm~\ref{alg:meshgen}.  The total number of refinements $M$ is a constant
to guarantee the termination of the algorithm.  Nevertheless, in
practice and in our implementation the refinement loop is terminated when a
quality criterion is satisfied, e.g., when the average minimum angle of the
triangles is above a predefined threshold.  The entire mesh generation process
for the spiral-shaped domain example is given in Figure~\ref{fig:spiralexample}.

\begin{algorithm}[H]
  \caption{Pseudocode for the triangular mesh generator}
  \label{alg:meshgen}
  \begin{algorithmic}[1]
    \Require{$\mathcal{C}$ is a sequence of corner points for a polygonal domain}
    \Require{$M$ is the total number of refinements}
    \Statex
    \Function{Adaptmesh}{$\mathcal{C}$}
    \Let{$\mathcal{T}_0^\prime$}{$\textsc{CDT}(\mathcal{C}$)}
    \Let{$\mathcal{T}_0$}{$\mathcal{T}_0^\prime$ with triangles outside of $\mathcal{C}$ removed}
    \For{$k \gets 1 \textrm{ to } M$}
    \Let{$\mathcal{T}_k^\prime$}{$\textsc{RGB}(\mathcal{T}_{k-1},\{\eta_T : T \in \mathcal{T}_{k-1}\})$}
    \Let{$\mathcal{T}_k$}{$\textsc{CPT}(\mathcal{T}_k^\prime)$}
    \EndFor
    \State \Return{$\mathcal{T}_N$}
    \EndFunction
  \end{algorithmic}
\end{algorithm}

\section{Implementation and example meshes}

We created a prototype of the mesh generator in Python for computational
experiments~\cite{adaptmesh2020}.  The implementation relies heavily on the
scientific Python ecosystem~\cite{virtanen2020scipy}.  It includes source code
from pre-existing Python packages \verb|tri| \cite{tri} (CDT implementation,
ported from Python 2) and the older MIT-licensed versions of \verb|optimesh|
\cite{optimesh} (CPT smoothing) and \verb|meshplex| \cite{meshplex} (edge
flipping).  Moreover, it performs adaptive mesh refinement using
\verb|scikit-fem| \cite{gustafsson2020scikit} and visualisation using
\verb|matplotlib| \cite{hunter2007matplotlib}.

\begin{figure}[htbp]
  \centering
  \includegraphics[width=0.32\textwidth]{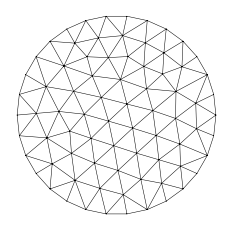}
  \includegraphics[width=0.32\textwidth]{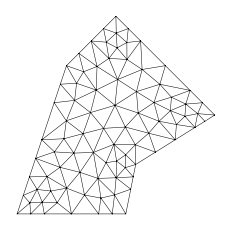}
  \includegraphics[width=0.32\textwidth]{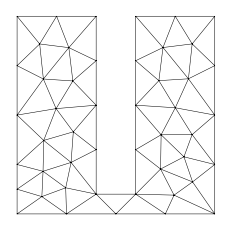}\\
  \includegraphics[width=0.22\textwidth]{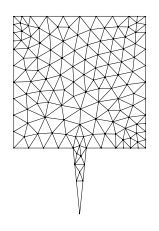}
  \includegraphics[width=0.32\textwidth]{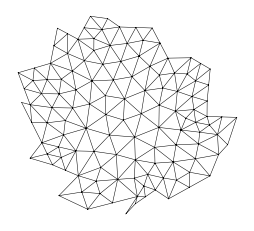}
  \caption{Some example meshes generated using our
    implementation of the proposed algorithm.
  The source code for the examples is available in \cite{tomgustafsson_2020}.}
\label{fig:moreexamples}
\end{figure}

Some example meshes are given in Figure~\ref{fig:moreexamples}.  By
default, our implementation uses the \emph{average} triangle quality\footnote{Triangle
quality is defined as two times the ratio of the incircle and circumcircle
radii.} above 0.9 as a stopping criterion which can lead to individual slit triangles.
This is visible especially in the last two examples that have small
interior angles on the boundary.

\section{An example application}

\begin{figure}
  \centering
  \includegraphics[width=0.45\textwidth]{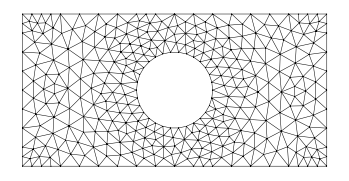}\\
  \includegraphics[width=0.5\textwidth]{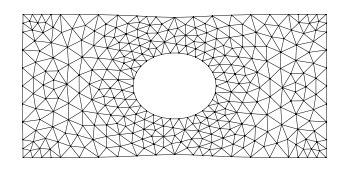}
  \caption{The original and the deformed meshes: the vertices of the original mesh are deformed using the
  finite element approximation of $\boldsymbol{u}$.}
  \label{fig:kirsch}
\end{figure}

This example utilises a variant of Algorithm~\ref{alg:meshgen} with steps 2 and 3
modified to allow for the inclusion of polygonal holes~\cite{adaptmesh2020}.
The holes are treated similarly as the sequence of corner points $\mathcal{C}$ in step 2
and the triangles inside the holes are removed in step 3.
We consider the domain
\[
\Omega = \{ (x, y) \in \mathbb{R}^2 : x^2 + y^2 \geq 1, -4 \leq x \leq 4,-2 \leq y \leq 2\}
\]
which is approximated by the triangular mesh given in Figure~\ref{fig:kirsch}.
We split the boundary of the domain into two as $\partial \Omega = \Gamma \cup (\partial \Omega \setminus \Gamma)$
where $\Gamma = \{ (x, y) \in \mathbb{R}^2 : |x|=4 \}$, and
consider the following linear elastic boundary value problem: find $\boldsymbol{u} : \Omega \rightarrow \mathbb{R}^2$ satisfying
\begin{alignat}{2}
  -\mathrm{div}\,\boldsymbol{\sigma}(\boldsymbol{u}) + \varepsilon \boldsymbol{u} &= \boldsymbol{0} \qquad &&\text{in $\Omega$},\\
  \boldsymbol{\sigma}(\boldsymbol{u})\boldsymbol{n} \cdot \boldsymbol{n} &= g \qquad &&\text{on $\Gamma$},\\
    \boldsymbol{\sigma}(\boldsymbol{u})\boldsymbol{n} \cdot \boldsymbol{t} &= 0 \qquad &&\text{on $\Gamma$},\\
  \boldsymbol{\sigma}(\boldsymbol{u})\boldsymbol{n} &= \boldsymbol{0} \qquad &&\text{on $\partial \Omega \setminus \Gamma$},
\end{alignat}
where
\[
\boldsymbol{\sigma}(\boldsymbol{u}) = 2 \mu \boldsymbol{\epsilon}(\boldsymbol{u}) + \lambda \,\mathrm{tr}\,\boldsymbol{\epsilon}(\boldsymbol{u}) \boldsymbol{I} \in \mathbb{R}^{2 \times 2}, \quad \boldsymbol{\epsilon}(\boldsymbol{u}) = \frac12\left(\nabla \boldsymbol{u} + \nabla \boldsymbol{u}^T\right) \in \mathbb{R}^{2 \times 2},
\]
and
$\boldsymbol{n} \in \mathbb{R}^2$ denotes the outward unit normal,
$\boldsymbol{t} \in \mathbb{R}^2$ is the corresponding unit tangent,
$\boldsymbol{I} \in \mathbb{R}^{2 \times 2}$ is the identity matrix,
and the remaining parameters are $\varepsilon = 10^{-6}$, $g = 10^{-1}$, $\mu = \lambda = 1$.

\begin{figure}
  \centering
  \includegraphics[width=0.5\textwidth]{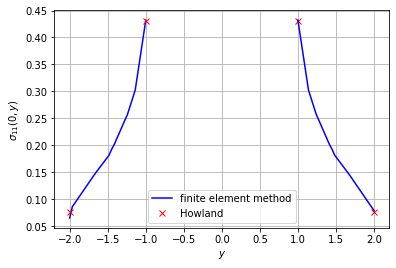}
  \caption{A comparison of the finite element approximation of the stress $\sigma_{11}(0, y)$ for $y \in [-2, 2] \setminus (-1, 1)$
  and the pointwise reference values from Howland~\cite{howland1930stresses}.}
  \label{fig:kirsch2}
\end{figure}

We solve the above problem using the finite element method and quadratic
Lagrange elements~\cite{Verf_rth_2013}.  The resulting deformed mesh, with the vertices of the original mesh displaced by the finite element approximation of $\boldsymbol{u}$,
is given in Figure~\ref{fig:kirsch}.
A reference value of the stress $(\boldsymbol{\sigma}(\boldsymbol{u}(x, y))_{11} = \sigma_{11}(x, y)$
at internal and external boundaries along $x=0$
is given by \mbox{Howland~\cite{howland1930stresses}}
via successive approximation. The reference values are $\sigma_{11}(0, 2) = \sigma_{11}(0, -2) \approx 0.75 g$ and
$\sigma_{11}(0, 1) = \sigma_{11}(0, -1) \approx 4.3 g$.
A comparison of the finite element approximation and the reference values
is given in Figure~\ref{fig:kirsch2}.

\section{Conclusions}
\label{sec:org615e973}

We introduced an algorithm for the generation of triangular meshes for explicit
polygonal domains based on the standard adaptive finite element method and
centroidal patch triangulation smoothing.  We presented a prototype
implementation which demonstrates that many of the resulting
meshes are reasonable and have an average triangle quality equal to or above
0.9.

A majority of the required components are likely to be available in an
existing implementation of the adaptive finite element method.  Therefore, the algorithm
can be a compelling candidate for supplementing an existing adaptive finite element solver 
with basic mesh generation capabilities.
Technically the approach extends to three dimensions although in practice the
increase in computational effort can be significant and the quality of the
resulting tetrahedralisations has not been investigated.

\bibliographystyle{siamplain}
\bibliography{mesh}

\end{document}